\documentclass[preprint]{elsarticle} %
\usepackage[english]{babel}
\usepackage{graphicx}
\usepackage{epsfig}
\usepackage{amssymb}
\usepackage{amsthm}
\usepackage{amsmath}
\usepackage{url}
\usepackage{pifont}
\newcommand{\C}{\mathbb C}

\def\lsim{\mathrel{\rlap{\lower4pt\hbox{\hskip1pt$\sim$}}
    \raise1pt\hbox{$<$}}}                

\newtheorem{corollary}{Corollary}

\newtheorem{proposition}{Proposition}
\newtheorem{theorem}{Theorem}

\input epsf.tex
\def\fig#1#2#3{
\vskip #3pc {\bf Fig.~#1.} {\it #2}  }
\def\psone#1#2{\centerline{
    \epsfxsize #2pc  \epsfbox{#1}
    }}
\def\pstwo#1#2#3#4{
   \centerline{
\epsfxsize=#2pc \epsfbox{#1}
      \qquad \epsfxsize=#4pc  \epsfbox{#3}
   }}%

\begin{document}
\begin{frontmatter}

\title{Solving polynomial eigenvalue problems by means of the Ehrlich-Aberth method}

\author[dm]{Dario A. Bini}

\ead{bini@dm.unipi.it}

\author[dm]{Vanni Noferini\corref{cor1}}
\ead{noferini@mail.dm.unipi.it}
\cortext[cor1]{Corresponding author.}
\address[dm]{Department of Mathematics, University of Pisa\\
Largo Bruno Pontecorvo 5, 56127 Pisa, Italy}

\begin{abstract}
Given the $n\times n$ matrix polynomial $P(x)=\sum_{i=0}^kP_i x^i$, we
consider the associated polynomial eigenvalue problem.
This problem, viewed in terms of computing the roots of the scalar
polynomial $\det P(x)$, is treated in polynomial form rather than in
matrix form by means of the Ehrlich-Aberth iteration.  The main
computational issues are discussed, namely, the choice of the starting
approximations needed to start the Ehrlich-Aberth iteration, the
computation of the Newton correction, the halting criterion, and
the treatment of eigenvalues at infinity. We arrive at an effective
implementation which provides more accurate approximations to the
eigenvalues with respect to the methods based on the QZ algorithm. The
case of polynomials having special structures, like palindromic,
Hamiltonian, symplectic, etc., where the eigenvalues have special symmetries
in the complex plane, is considered. A general way to adapt the
Ehrlich-Aberth iteration to structured matrix polynomial is
introduced. Numerical experiments which confirm the effectiveness of
this approach are reported.

\bigskip


\noindent
\end{abstract}


\end{frontmatter}

\section{Introduction}
Given two positive integers $k, \ n$ and matrices $P_j\in\C^{n\times n}$,
$j=0,\ldots,k$ consider the matrix polynomial

\begin{equation}\label{pol}
P(x)=\sum_{j=0}^k P_j x^j
\end{equation}
where $P_k\ne 0$, so that $P(x)$ has degree $k$, and define the scalar
polynomial $p(x)=\det P(x)$ of degree $N\le nk$. Assume that $P(x)$ is
regular, that is $p(x)$ is not identically zero.

Under such an assumption, the polynomial eigenvalue problem associated with $P(x)$ consists in
computing the roots of the polynomial $p(x)$ which are called the
eigenvalues of the matrix polynomial $P(x)$. Observe that, if $P_k$
has not full rank, then $N<nk$. In this case, it is convenient to
introduce $nk-N$ eigenvalues at infinity and say that the matrix
polynomial $P(x)$ has $nk$ eigenvalues including the $nk-N$
eigenvalues at infinity.

Our interest is addressed to the design and analysis of efficient
algorithms for the polynomial eigenvalue problem based on the
Ehrlich-Aberth iteration \cite{aberth, borschsupan, ehrlich}.

Recently, much literature has been addressed to the polynomial
eigenvalue problem (PEP). For the numerical solution 
of PEPs, fast and numerically stable methods are sought.  Several
algorithms have been introduced based on the technique of
linearization where the polynomial problem is replaced by a linear
pencil with larger size and the customary methods for the generalised
eigenvalue problem are applied. For more details,
see for instance \cite{bestqz, MMMM, vecspa, mt}
and the references therein.

Specific attention concerns structured problems where the matrix
coefficients $P_j$ have some additional property which is reflected on
structural properties of the roots. For instance, in the case of
T-palindromic polynomials  \cite{lanccnal, urv}, 
where $P_j=P_{k-j}^T\in\C^{n\times
  n}$, the roots are encountered in pairs $(x,1/x)$. In general, we
may consider the case of structures where the roots can be grouped in
pairs as $(x,f(x))$, where $f(x)$ is any analytic function such that
$f(x)=f^{-1}(x)$ \cite{nuovo}. In this case the goal is to design
algorithms which take advantage of this additional information about
the eigenvalues and deliver approximations to the eigenvalues which
respect these symmetries independently of the rounding
errors.

The Ehrlich-Aberth iteration (EAI) was historically first mentioned in \cite{borschsupan} and afterwards
 independently rediscovered many times. It is one of the many simultaneous
iteration techniques available in the literature for the numerical
approximation of polynomial roots \cite{macnamee, petkovic}. In \cite{B,BF} the
EAI has been combined with various techniques like the Rouch\'e
theorem, the Newton polygon technique, and the Gerschgorin inclusion
theorems for arriving at efficient and robust software
implementations. The package Polzeros, designed in \cite{B}, provides
a robust and reliable tool for approximating roots of polynomial in
floating point arithmetic. The package MPSolve designed in \cite{BF}
provides certified approximations to any desired number of digits of
the roots of any polynomial.

The EAI has been used in \cite{bgt} to solve the generalised
tridiagonal eigenvalue problem where the software provides effective
accelerations in terms of CPU time. It has been used in
\cite{plestenjak} for quadratic hyperbolic tridiagonal eigenvalue
problems.

In this paper we present an adaptation of the Ehrlich-Aberth method
for the numerical solution of PEPs.  The
main computational issues that we analyze are the choice of the
starting approximations, the computation of the Newton correction, the
halting criterion, the design of {\em a posteriori} error bounds, and
the management of the problematic (multiple) eigenvalues at zero and
at infinity.

Concerning the choice of the starting approximations, we propose a
generalization to matrix polynomials of a technique introduced in
\cite{B} for scalar polynomials based on the Rouch\'e theorem. In
fact, we rely on the generalization to matrix polynomials of the
Rouch\'e theorem given in \cite{ma} and provide a way to determine an
annulus in the complex plane which contains all the eigenvalues.

The Newton correction is computed by means of Jacobi's formula for the
differential of the determinant of any square matrix in terms of the
trace of the matrix $P(z)^{-1}P'(z)$.  The halting condition is given
in terms of the condition number of $P(z)$. Eigenvalues at zero and at
infinity can be removed, to a certain extent, by using the specific
features of the EAI relying on the information provided by the
singular values of $P_0$ and $P_k$. A posteriori error bounds are
given by constructing a set of inclusion disks relying on the
Gerschgorin theorem and adapting the results of \cite{smith} to the
case of matrix polynomials.

The computational analysis of this method shows that the number of
arithmetic operations (ops) is $O(k^2n^3+kn^4)$. In the case where the
degree $k$ is large with respect to the square root of the matrix size $n$, this
complexity bound compares favourably with the bound $O(k^3n^3)$ of the
customary matrix-based algorithms like the QZ applied to a linearization.
Cases of this kind can be encountered, for instance,
 in the truncation of matrix power series \cite{wilson}.

The EAI does not compute the eigenvectors, which are sometimes needed as well. 
Nevertheless, after a good 
approximation of the eigenvalues is obtained, other methods, 
e.g. the SVD or the inverse iteration, can be used to compute the 
eigenvectors without increasing the complexity of the algorithm. We were 
able to compute eigenvectors with high accuracy using the eigenvalues 
given by the EAI.

We consider the case of polynomials endowed with specific properties
like palindromic, T-palindromic, Hamiltonian, symplectic polynomials,
whose eigenvalues have special symmetries in the complex plane. We
propose a unifying treatment of this class of structured polynomials
and show how the EAI can be adapted to deal with these classes in a
very effective way. In fact, our variant of the EAI enables one to
compute only a subset of eigenvalues and to recover the remaining part
of the spectrum by means of the symmetries satisfied by the
eigenvalues. By exploiting the structure of the problem, this
approach leads to a saving on the number of operations and provides
algorithms which yield numerical approximations fulfilling
the symmetry properties.

We conclude our discussion by presenting the results of several
numerical experiments performed in order to test the effectiveness of
our approach in terms of speed and of accuracy.  We have compared the
Ehrlich-Aberth iteration with the Matlab functions {\tt polyeig} and
{\tt quadeig}\cite{quadeig}. In the structured case, we have also considered, when
available, other structured methods, say, the URV algorithm by
Schr\"{o}der \cite{urv}.

We show that the EAI is much faster than the available techniques in
the case where the degree is larger with respect to the size of
the matrices. Moreover, for the test problems NLEVP of \cite{nlevp},
it turns out that the accuracy of the computed approximations is
generally better than the accuracy obtained with the available
algorithms.

The paper is organised as follows. In Section \ref{sec2} we recall the
Ehrlich-Aberth method and discuss the main computational issues
encountered in its implementation.
In Section \ref{sec:struct} we consider the case of ``structured polynomials'', i.e. 
the case where the matrix coefficients have some special
properties. 
Section \ref{numex} reports the results of the numerical experiments.

\section{The Ehrlich-Aberth method for matrix polynomials}\label{sec2}

Given a vector $y^{(0)}\in\C^N$ of initial guesses for the $N$ roots
of the polynomial $p(x)$, the EAI provides the sequence of
simultaneous approximations $y^{(i)}$ given by
\begin{equation}\label{eaj}
\begin{array}{l}
y^{(i+1)}_j=y^{(i)}_j- \displaystyle \frac{\mathcal N(y_j^{(i)})}
{1-\mathcal  N(y_j^{(i)})\mathcal A_j(y^{(i)})}, ~~~\mathcal
  N(x)=\frac{p(x)}{p'(x)}\\[3ex]
\mathcal A_j(y^{(i)})=\displaystyle 
\sum_{\ell=1,\ell\ne j}^N\frac1{y_j^{(i)}-y_\ell^{(i)}}
\end{array}
\end{equation}
where $\mathcal N(x)$ is the Newton correction. It is easy to check
that the $j$th update in \eqref{eaj} is nothing else but the Newton
iteration applied to the rational function
$p(x)/\prod_{\ell=1,~\ell\ne j}^N(x-y_{\ell}^{(i)})$, so that the EAI
provides a way to implement the implicit deflation of the roots.

Besides the Jacobi-style version of EAI we may formulate the
Gauss-Seidel version of EAI, that is
\begin{equation}
\begin{array}{l}\label{eags}
y_j^{(i+1)} = y_j^{(i)} -\displaystyle \frac{\mathcal N(p(y_j^{(i)}))}
{1-\mathcal N(p(y_j^{(i)})) \mathcal A_j( y^{(i)},  y^{(i+1)}) },\\[2ex]
\displaystyle \mathcal A_j( y^{(i)},  y^{(i+1)})=
\sum_{\ell=1}^{j-1}\frac{1}{y_j^{(i)}-y_{\ell}^{(i+1)}} +
\sum_{\ell=j+1}^{N}\frac{1}{y_j^{(i)}-y_{\ell}^{(i)}}.
\end{array}
\end{equation}
The method, in the Jacobi version, is known to converge cubically for
simple roots and linearly for multiple roots \cite{petkovic}. In the
Gauss-Seidel version, convergence is slightly faster.  In practice,
good global convergence properties are observed; a theoretical
analysis of global convergence, though, is still missing and
constitutes an open problem.

With the term {\em vector iteration} of the Ehrlich-Aberth method we
refer to the step which provides the vector $y^{(i+1)}$ given the
vector $y^{(i)}$.  We use the term {\em scalar iteration} for
indicating the single step performed on the generic scalar component
of the vector $y^{(i)}$.
 
In the case of a scalar polynomial of degree $N$ the cost of a scalar
iteration is $O(N)$ arithmetic operations. In this way, the cost
of a vector iteration is at most $O(N^2)$ ops and is substantially
reduced when most components have been numerically approximated, so
that few scalar iterations must be performed in order to carry out the
vector iteration.

The number of scalar iterations needed by the floating point
implementation in order to find approximations which are exact roots
of a slightly perturbed polynomial is, in practice, $O(N)$ if the
starting approximations are computed by means of the Newton polygon
technique \cite{B}. This technique is particularly effective when the
polynomial has roots with moduli which are very unbalanced.

Crucial aspects for an effective implementation of the EAI to matrix
polynomials are
\begin{enumerate}
\item the computation of the Newton correction $p(x)/p'(x)$ given the
  value of $x$ and of the input coefficients $P_j$, $j=0,\ldots,k$;
\item a criterion for stopping the iterations;
\item the choice of the initial approximations.
\end{enumerate}

\subsection{Computing the Newton correction}

In the literature, methods based on some factorizations of $P(x)$ 
were developed to compute the Newton correction for functions 
that have the same zeros of $p(x)$: 
e.g., the method in \cite{kublanovskaya}, later proved to lack 
theoretical rigour and corrected in \cite{kainshingal}. Other kinds of Newton-like 
approaches were presented in \cite{jain2, ruhe}.

If one wishes to work with $p(x)$ itself, a naive way to compute the 
Newton correction $p(x)/p'(x)$ would be to evaluate first the coefficients of the
polynomial $p(x)$, say, by means of the evaluation-interpolation
technique, and then to apply right after the Ehrlich-Aberth method to
the scalar equation $p(x)=0$. This approach would however come across
numerical problems due to numerical instability and to overflow and
underflow situations encountered in the computation of determinants.

It is therefore wise to conceive a strategy to avoid the explicit
calculation of the coefficients of $p(x)$. 

An effective way rests upon
the well-known Jacobi's formula for the differential of the
determinant of any invertible square matrix $A$:

\begin{equation}\label{jacobiformula}
d (\ln (\det A)) = \mathrm{tr}(A^{-1} d A),
\end{equation}
where tr denotes the trace.  This way, we obtain the following
expression for the derivative $p'(x)$

\[
p'(x)=\frac{d \det (P(x))}{d x} = 
\det(P(x)) \cdot \mathrm{tr}(P(x)^{-1}\cdot \frac{d P(x)}{d x}).
\]

This formula allows us to evaluate the Newton correction $p(x)/p'(x)$,
which is the centerpiece for the EAI, without explicitly calculating
$p(x)$:

\begin{equation}\label{newtcor}
p(x)/p'(x)=\frac 1{\mathrm{tr}(P(x)^{-1}P'(x))}.
\end{equation}

An evaluation of $P(x)$ and $P'(x)$ by means of Horner's method,
followed by a numerical matrix inversion, allows to compute the trace
of $P(x)^{-1}P'(x)$ in $\mathcal{O}(k n^2 + n^3)$ operations.

We mention that, even though we independently formulated it, we found out later on that the possible use 
of \eqref{jacobiformula} in a numerical method for PEPs 
had been already suggested in \cite{WG, lanlam}. However, in these instances it was proposed to use 
it to apply the Newton method to approximate each single eigenvalue in sequence, mentioning the
possibility to use an implicit deflation of 
previously found roots \cite{maehly}, in order to avoid that the method converges twice to the same eigenvalue. 
This leads to a formula akin to 
\eqref{eaj}, \eqref{eags}, with the difference that the summation in the term $\mathcal A$ is performed only up to 
the number of roots that have already been approximated. 
This is a crucial detail, because such a sequential implementation of the Newton method does not seem to achieve the 
same efficiency with respect to the EAI.

\subsection{Stopping criterion}
At the generic $i$th vector iteration it is crucial to decide whether the update
of the $j$th component of the vector $y^{(i+1)}$ must be performed or
the scalar iteration in that component must be halted.

Observe that, if $\xi$ is a root of $p(x)$, that is $\det P(\xi)=0$,
then as the approximation $x$ gets close to $\xi$ the matrix $P(x)$
becomes ill-conditioned. This makes quite natural to stop the
iterations if the reciprocal of the condition number $\mu(P(x))$ is
less than a prescribed tolerance $\tau_1$. This criterion makes sense if the eigenvalue that we want
to approximate is semi-simple. In the case of a defective eigenvalue
$\lambda$ with Jordan chains of length at most $m$, in view of the
results in \cite{pp}, it is more convenient to stop the iterations if
the reciprocal of $\mu(P(x))$ is less than $\tau_1^m$. This latter
condition is hard to implement since it is not easy to evaluate
numerically the length of the Jordan chains of a matrix polynomial.
Using the former stopping criterion may lead to a premature halt of the
algorithm in the case of defective eigenvalues.

As an alternative to the previous stopping condition, following \cite{tisseurduemila}, define $\alpha(x) = \sum_{\ell=0}^k |x^{\ell}|$. 
If $y_j^{(i)}$ is not an eigenvalue of $P(x)$ then the quantity
\[
\eta(y_j^{(i)})=\left( \parallel \left( P(y_j^{(i)}) \right)^{-1} \parallel_2 (1+ \alpha(y_j^{(i)})) \right)^{-1}
\]
measures the backward error for the approximation $y_j^{(i)}$, and can be cheaply evaluated during the EAI. The iteration can be 
halted when $\eta(y_j^{(i)})$ is smaller than a given tolerance.

 For simple  eigenvalues, no significant differences emerged between the two alternative 
possibilities. Therefore, our default choice was in favour of the criterion based on the 
condition number. 

It is also convenient to add, with the ``or'' logic operator, the
following condition
\begin{equation}\label{stop2}
|\mathcal N(y_j^{(i)})/(1-\mathcal N(y_j^{(i)})\mathcal A_j
(y^{(i)},y^{(i+1)})|\le\tau_2|y_j^{(i)}|
\end{equation}
where $\tau_2$ is a given tolerance. This condition says that
the computed correction is too tiny and would not change the significant digits of
the current approximation.

\subsection{Choosing initial approximations}

As pointed out in \cite{aberth, B, gug}, practically
effective choices of initial approximations for the EAI are complex
numbers equally displaced along circles. For instance, in
\cite{aberth} it is proposed to choose initial approximations
displaced along a circle centered at the origin of sufficiently large
radius so that it contains all the roots.  In \cite{gug} the radius of
the circle is suitably chosen.  This strategy does not work
effectively for polynomials having zeros with very large and with very
small moduli.  In \cite{B} this drawback is overcome by considering
different circles centered at the origin of suitable radii. The
computation of these radii relies on the Rouch\'e theorem.

Here we try to extend this technique to a certain extent.  We recall
that, according to the Rouch\'e theorem, if $s(x)$ and $q(x)$ are two
polynomials such that
\[
|s(x)|>|q(x)|,~~\hbox{for}~~|x|=r,
\]
then $s(x)$ and $s(x)+q(x)$ have the same number of roots in the open
disk $\{z\in\C:\quad |x|<r\}$. Applying this property with $s(x)=x^m$
and $q(x)=p(x)-s(x)$, for $0\le m\le N$, implies that if
$r^m>\sum_{j=0,j\ne m}^N |a_j|r^j$ then the polynomial $p(x)$ has $m$
roots in the open disk of center $0$ and radius $r$. This property is
at the basis of the criterion described in \cite{B}, based on the
Newton polygon construction, for choosing initial approximations
equidistributed along different circles centered in $0$.

In order to extend this criterion to the case of matrix polynomials we
need a generalisation of the Rouch\'e theorem to matrix
polynomials. We report the following result of \cite{ma} which we
rephrase in a simpler way better suited for our problem.

\begin{theorem} Let $S(x)$ and $Q(x)$ be matrix polynomials and let $r$ 
be a positive real. If $S(x)^*S(x)-Q(x)^*Q(x)$ is positive definite
for $|x|=r$, then the polynomials $\det S(x)$ and $\det(S(x)+Q(x))$
have the same number of roots of modulus less than $r$.
\end{theorem}

The following result is an immediate consequence of the above theorem
applied to the polynomial $P(x)$ of \eqref{pol} with $S(x)=x^m P_m$
and $Q(x)=\sum_{i=0,\, i\ne m}^{k}x^iP_i$.

\begin{corollary}\label{cor}
Assume that 
\begin{equation}\label{succeq}
P_m^*P_m r^{2m}-(\sum_{j=0,j\ne m}^k P_j^*\bar
x^j)(\sum_{j=0,j\ne m}^k P_j x^j)\succ 0,~~\hbox{for}~~|x|=r,
\end{equation}
where $A\succ B$ means that $A-B$ is positive definite.
 Then the matrix polynomial $P(x)$ has $mk$ eigenvalues in the
open disk of center $0$ and radius $r$.
\end{corollary}

Observe that if $\det P_m=0$ then condition \eqref{succeq} cannot be
verified. In fact, the vector $v$ such that $P_mv=0$ would be such
that
\[
v^*(\sum_{j=0,j\ne m}^k P_j^*\bar x^j)(\sum_{j=0,j\ne m}^kP_j
x^j)v\le 0
\]
 which is absurd.

In particular, if $\det P_k\ne 0$ the above corollary, applied with
$m=k$, implies that all the eigenvalues of $P(z)$ are included in the
disk of center $0$ and radius $r$ provided that
\begin{equation}\label{r}
r^{2k}
P_k^*P_k-(\sum_{j=0}^{k-1}P_j^* \bar {x}^j)(\sum_{j=0}^{k-1}P_j x^j)\succ 0,\quad\hbox{for }|x|=r.
\end{equation}
Observe that the latter condition is implied by
\begin{equation}\label{rr}
r^{2k} P_k^*P_k\succ \sum_{j=0}^{k-1}r^jP_j^*P_j+I\sum_{j>
i}\rho(|P_j|^*|P_i|+|P_i|^*|P_j|).
\end{equation}

Similarly, applying Corollary \ref{cor} with $m=0$ provides a disk
where $P(x)$ has no eigenvalues.

As an example of application, consider the $5\times 5$ quadratic
matrix polynomial $P(x)=Ax^2+Bx+A^T$, where $B$ is the tridiagonal
matrix defined by the entries $[1,2,1]$, and $A$ is the matrix with
diagonal entries $100, 1, 1/1000, 1/100000$, superdiagonal entries
equal to 1 and with zero entries elsewhere.  The eigenvalues of $P(x)$
have moduli {\tt 2.0050e+05, 1.4969e+03, 1.0000e+00, 1.0000e+00,
  1.0000e+00, 1.0000e+00, 6.6805e-04, 4.9874e-06}. The criterion based
on the above corollary in the form \eqref{succeq} yields the bound
${\tt 4.4e-6}<|x|<{\tt 2.24e5}$ which is quite good.  Applying
condition \eqref{rr} yields the bounds ${\tt 1.96e-6}<|x|<{\tt 5.1e5}$
which is still good.

Similar results have been obtained in \cite{sharify} in the framework
of tropical algebras.  A different heuristic approach, which relies on Corollary
\ref{cor}, is to select the values of the radii by 
considering the
inequality $\|P_m\|_2r^m>\sum_{i=0,\, i\ne m}^{k}\|P_i\|_2r^i$  in place of \eqref{succeq}.
  This strategy, applied in the form
$r^k>\sum_{i=0}^{k-1}\|P_i\|_2r^i$
 leads to the  criterion based on computing the Newton
polygon of the polynomial $\sum_{i=0}^k \|P_i\|_2x^i$. This is
the default choice of the starting approximations performed in our
implementations.

\subsection{A posteriori error bounds}
In the case of a scalar polynomial $p(x)$ of degree $N$, given a set of
approximations $x_1,\ldots ,x_N$ to the roots of $p(x)$ it is possible
to prove that \cite{smith} the set of disks $D_i=D(x_i,r_i)$ of
center $x_i$ and radius $r_i=n|p(x_i)/(p_N \prod_{j=1,\, j\ne i}^N
(x_i-x_j))|$ is such that
\begin{enumerate}
\item the union of the disks contains all the roots of $p(x)$
\item each connected component formed by the union of, say, $c$
  overlapping disks, contains $c$ roots of $p(x)$.
\end{enumerate}
The set formed by $D_i$, $i=1,\ldots,N$ with the above properties is
said {\em set of inclusion disks}.

In the case of a matrix polynomial $P(x)$ where $\det P_k\ne 0$, it is
quite cheap to compute a set of inclusion disks. In fact, if $P(x)=\Pi
L U$ is the PLU factorization of $P(x)$, then $|p(x)|=|\det
P(x)|=\prod_{j=1}^n |u_{j,j}|$, where $U=(u_{i,j})$.  Moreover, the
leading coefficient $p_N$ of $\det P(x)$ coincides with $\det P_k$
which can be computed once for all.  Observe that the LAPACK routine
{\tt zgesv} which solves a linear system with the matrix $P(x)$, used
to compute the Newton correction $1/\hbox{trace}(P(x)^{-1}P'(x))$,
applied with $x=x_i$, provides at a negligible cost also the radius
$r_i$.

The availability of a set of inclusion disks enables one to perform a cluster analysis. 
In fact, once an isolated disk has been detected, we have isolated a single eigenvalue of the matrix polynomial $P(x)$. Once we have 
detected a set of $c$ overlapping disks isolated from the remaining inclusion disks, we have detected a cluster formed by $c$ 
eigenvalues of $P(x)$. 

A different {\em a posteriori} error bound can be obtained by using a
classical result \cite{henrici}. The disk of center $x_i$ and
radius $\hat r_i=n|p(x_i)/p'(x_i)|$ contains a root of the polynomial
$p(x)$. However, the set of disks obtained in this way does not
fulfill properties 1 and 2 of the set of inclusion disks. It is worth
pointing out that the computation of $\hat r_i$ is inexpensive since
the Newton correction $p(x_i)/p'(x_i)$ is computed by the EAI.
Moreover, this {\em a posteriori} error bound still holds if the
leading coefficient $P_k$ is singular.

\subsection{Multiple eigenvalues}
Computational difficulties may be encountered in the case of multiple
eigenvalues. In fact, the rate of convergence for multiple eigenvalues
is linear, with respect to the cubic behaviour for simple eigenvalues.
Moreover, for defective eigenvalues the standard stop condition may
lead to a premature halt. For this reason, if it is possible to detect \emph{a priori} 
multiple eigenvalues, it is advisable to deflate them; if it is not possible to spot all of 
them theoretically, even lower bounds on the multiplicity are very helpful. If multiple 
eigenvalues are not predicted theoretically, one must rely on the cluster 
analysis to identify them and modify accordingly the stopping criterion.

A common situation that leads to multiple eigenvalues is met when the
extremal coefficients are rank-deficient matrices. In this case, $0$
and/or $\infty$ have multiplicity greater than or equal to $1$. This situation can
be  circumvented to a certain extent. 

In the case of $m$ eigenvalues at infinity, 
one may just start with an approximation vector $y$ of length $nk-m$, 
acknowledging that the
determinant $p(x)$ has in fact degree $nk-m$; if there are $m$ zero eigenvalues
it is possible to set to zero $m$ components of the vector $y^{(0)}$
avoiding to update them. 

The number of null singular values of $P_0$ provides a lower bound to
the number of null eigenvalues of $P(x)$. Similarly, the number of zero 
singular values of $P_k$ provides a bound to the number of eigenvalues
at infinity. This way, the precomputation of the SVD of $P_0$ and $P_k$ may 
increase the performance of the EAI. Equivalently, one may perform any rank-revealing 
factorization (e.g., QR) instead of the SVD. Sometimes, the structure of the coefficients allows 
to achieve better bounds 
(e.g., if the same rows/columns of many consecutive extremal coefficients are zero).

In our implementation, the rank of the extremal coefficients is tested. If it is less than $n$, it is also checked if $P_0$ and $P_1$ 
(resp., $P_{k}$ and 
$P_{k-1}$) share any common zero row/column. Thus, any manifest presence of zero and infinite eigenvalues is 
exploited, forcing deflation of all the guaranteed roots.
Moreover, if the test detects the presence of eigenvalues at $0$ (resp., 
$\infty$), in order to avoid a premature stop for other undetected eigenvalues at $0$ (resp., $\infty$) 
the stopping criterion is made stricter. The stronger stop condition 
requires that, for eigenvalues smaller (resp., larger) than a given bound, 
either the relative correction criterion \eqref{stop2} is satisfied with tolerance $\tau_2$ or the relative correction criterion \eqref{stop2} is satisfied with
tolerance $\tau_2^{1/2}$ and, simultaneously, the reciprocal condition number criterion is satisfied with tolerance $\tau_1$. This 
heuristic device worked very effectively in our experiments, leading to satisfying results also in problems with multiple eigenvalues 
at either zero or infinity (see Section \ref{numex}). 

If the leading coefficient $P_k$ is singular and if the degree of
$p(x)=\det P(x)$ is not available together with the leading
coefficients of $p(x)$, then it is not possible to generate a set of
inclusion disks and to perform a cluster analysis. However, in this
case we may apply an effective technique based on a rational
transformation of the variable $x$. For instance, the variable $x$ is
replaced by the M\"obius function $x=x(z)=(\alpha z+\beta)/(\gamma
z+\delta)$ such that $\alpha\delta-\gamma\beta\neq 0$, and the
polynomial $P(x)$ is replaced by the polynomial $Q(z)=(\gamma
z+\delta)^kP(x(z))$. This way the eigenvalues at infinity of $P(x)$
are mapped into eigenvalues of $Q(z)$ at $-\delta/\gamma$. Moreover,
$Q(z)$ has no eigenvalues at infinity provided that $\alpha/\gamma$ is
not eigenvalue of $P(x)$.  The substitution of variable can be
performed implicitly without actually computing the coefficients of
$Q(z)$ except for $Q_k$. We refer the reader to Section \ref{sec:struct}
and to \cite{vanni} for more details.

\subsection{Linearization as a possibility}\label{linear}
In the paper \cite{dicksonpal} the computation of the Newton
correction was carried out by first linearizing the polynomial by
means of companion-like pencil, and then by evaluating the trace by
means of an LQ factorization.  This approach has a computational
complexity of $\mathcal{O}(n^3+kn^2)$ ops per scalar iteration, which
is the same that is achieved by our algorithm.  However, this cost can
be reduced to $O(n^2k^2)$ ops if the matrix pencil is reduced to
triangular-Hessenberg form before the Jacobi formula is applied.  This
fact is the main advantage of using a linearization.  However,
linearization techniques, if not properly used, may lead to an
undesired increasing of the eigenvalue condition numbers
\cite{bestqz}.

\section{The case of structured polynomials}\label{sec:struct}
The EAI is particularly suited to deal with matrix polynomials endowed
with specific structures of the matrix coefficients. We are interested
in matrix structures which induce particular symmetries on the
location of the eigenvalues.  Polynomials of this kind are encountered
in the applications and include, for instance, palindromic,
T-palindromic, symplectic and Hamiltonian polynomials.

Customary PEP-solving algorithms, such as the application of the QZ to
any suitable linearization of the polynomial, are not able to fully
catch these symmetries of the spectrum. In the literature, there are
specific matrix methods that achieve this goal.  The EAI enables to
exploit the additional information both in the computation of the
Newton correction and in the choice and in the management of the
(initial) approximation of the roots in view of the structure-induced
symmetries.  We will see this later on.

Assume that the structured PEP is such that the eigenvalues appear in
pairs $\{x,f(x)\}$, with $f(f(x))=x \ \forall x$. A naive adaptation of
  the EAI to this property would be to apply \eqref{eaj} or
  \eqref{eags} updating only the first half of the components of the
  vector $y$ and simultaneously imposing $y^{(i)} = f(y^{(i-nk/2)})$,
  $i=nk/2+1,\dots,nk$. In numerical experiments, this approach seems not
  to be always efficient in terms of number of scalar iterations needed for
  numerical convergence. This motivates the design of more
  sophisticated structured variants of the EAI, that we are going to describe in the 
following.

Before analyzing the various classes of structured matrix polynomials,
we recall some basic definitions of special matrices.

An $n\times n$ square matrix $A$ is said to be \emph{symmetric} if
$A^T=A$ and \emph{skew-symmetric} if $A^T=-A$.  Let $n=2m$. The matrix
$A$ is said to be \emph{Hamiltonian} if it is such that $A^T J =-J A$
where $J$ is the matrix $\left(\begin{smallmatrix} 0 & I_m\\ -I_m &
  0\\
\end{smallmatrix}\right)$;
$A$ is said to be
  \emph{skew-Hamiltonian} if it is such that $A^T J = J A$;
$A$  is said to be
  \emph{symplectic} if it is such that $A^T J A =J$. 

Every skew-Hamiltonian matrix can be obtained as the square of a
Hamiltonian matrix, and conversely the square of a Hamiltonian matrix
is always skew-Hamiltonian\cite{fmm}. Symplectic matrices are exponential of
Hamiltonian matrices.

\subsection{Skew-Hamiltonian and even-dimensional skew-symmetric}
A skew-symmetric polynomial $P(x)$ is a polynomial whose coefficients $P_j$,
for $j=0,\dots,k$,
are skew-symmetric constant matrices. If the
coefficients have even size $n=2m$, we say that $P(x)$ is an even-dimensional
skew-symmetric polynomial. A skew-Hamiltonian polynomial is defined as
a polynomial whose coefficients $P_j$ are all skew-Hamiltonian
matrices. Classical eigenvalue problems for skew-Hamiltonian matrices \cite{vanloan} are a 
special case of skew-Hamiltonian PEPs.

These two classes of polynomials are closely related, because
multiplication by $J$ maps one class onto the other. A common feature
is that the spectrum of any polynomial in these two classes contains
only eigenvalues of even multiplicity. In fact, 
the determinant of a
matrix polynomial $P(x)$ belonging to these two classes can be written
as
\[
p(x)=\det P(x)=q(x) \cdot q(x),
\] 
for a suitable polynomial $q(x)$. For the special case of a real
skew-symmetric matrix pencil, a proof was  given in \cite{laro}
where a special Kronecker form was derived. The more general case comes from classical 
results on determinants \cite{muir}. Let us give here a simple proof of the statement for an
even-dimensional skew-symmetric complex matrix polynomial using modern terminology.

\begin{proposition}
Let $P(x)=-P(x)^T$ be a $2m \times 2m$ skew-symmetric matrix
polynomial. Then $p(x)=\det P(x) = q(x) \cdot q(x)$ for some scalar
polynomial $q(x)$.
\end{proposition}
\emph{Proof.} We shall prove the proposition by induction on $m$. For
$m=1$ the statement is obvious.
Suppose now that any $(2m-2) \times (2m-2)$ skew-symmetric polynomial
has the desired property. Let $\Pi$ be a $2m \times 2m$ permutation
matrix and let $Q(x):=\Pi P(x) \Pi^T$. Suppose that $\Pi$ is such that
\[
Q_0(x):=Q(1:2,1:2)=:\left(\begin{array}{cc} 0 & r(x)\\ -r(x) & 0\\
\end{array}\right)
\]
 is nonsingular, where $r(x)$ is a suitable 
nonzero scalar polynomial. Notice that such an assumption can be
safely made because if that was false for any $\Pi$ then $P(x)=0$ 
so $p(x)=0$ and there would be nothing to prove.
Now let $ Q(x)=\left(\begin{smallmatrix} Q_0(x) & A(x)\\ -A(x)^T &
  Q_1(x)\\
\end{smallmatrix}\right)
$, 
where the polynomial matrices $A(x)$ and $Q_1(x)$ have, respectively,
dimensions $2 \times (2m-2)$ and $(2m-2) \times (2m-2)$; also, let
$\rho(x):=r(x)^{m-1}$. Define the rational function
$S(x):=Q_1(x)+\frac{1}{r(x)} A(x)^T\left(\begin{smallmatrix} 0 &
  1\\ -1 & 0\\
\end{smallmatrix}\right)A(x).$ Clearly, $r(x) S(x)$ is a $(2m-2) \times (2m-2)$
skew-symmetric matrix polynomial; therefore, by the inductive hypothesis, $\det
S(x)=\frac{\theta(x)^2}{\rho(x)^2}$ where $\theta(x)$ is a suitable scalar
polynomial.
Moreover, $S(x)$ is the Schur complement of $Q_0(x)$. Thus,
$p(x)=\frac{r(x)^2 \theta(x)^2}{\rho(x)^2}$, so $p(x)$ is the square
of some scalar rational function
$q(x)=\frac{\theta(x)}{r(x)^m}$. Since $p(x)$ is a polynomial, $q(x)$
must be a polynomial as well. $\Box$\\

This is a particularly useful property which can be fully exploited by
the Ehrlich-Aberth method. In fact, instead of applying the EAI to the
polynomial $p(x)$ of degree $2mk$, one can apply the EAI to the
polynomial $q(x)$ of degree $mk$ even though $q(x)$ is not explicitly
known.

More precisely, since $p'(x)/p(x) = 2 q'(x)/q(x)$,
one can compute the Newton correction $q(x)/q'(x)$ by means of
\[
q(x)/q'(x)=2p(x)/p'(x)=2/\hbox{tr}(P(x)^{-1}P'(x)).
\]

This way, the length of the vector of the approximations $y$ in \eqref{eaj} 
or in \eqref{eags} is reduced
from $2mk$ to $mk$, moreover, the skew-Hamiltonian or the
skew-symmetric structure of the coefficients can be exploited in the
computation of $P(x)^{-1}P'(x)$.

\subsection{Palindromic and symplectic}
The polynomial $P(x)$ is called \emph{purely palindromic} if
$\mathrm{Rev}P(x)=P(x)$, where the reversal polynomial
$\mathrm{Rev}(P(x))$ is defined by $\mathrm{Rev}P(x):=x^k P(x^{-1})$.
The polynomial $P(x)$ is called \emph{T-palindromic} if
$\mathrm{Rev}P(x)=P(x)^T$. Both these structures induce a symmetry
$(x,1/x)$ in the spectrum. There is a vast literature on this kind of
structure; see, e.g., \cite{lanccnal, hybrid, MMMM, urv} and 
the references therein. The same structure appears in
the standard eigenvalue problem for a symplectic matrix\cite{fass}.

For this class of PEPs the change of
variable $z:=x+1/x$ is useful. In \cite{dicksonpal} it was shown that
the use of a non-standard polynomial basis, called the Dickson basis,
leads to a suitable linearization of the purely palindromic
case. Moreover, it was shown that if $P(x)$ is T-palindromic then it
is possible to build a new skew-Hamiltonian matrix polynomial $M(z)$
such that $\det M(z(x))=p(x) \cdot p(x)$; the Dickson basis was then
used to obtain a useful linearization. In the following, we will show 
how to avoid the explicit use of the Dickson transformation.  This has
the advantage that there is no potential loss of accuracy for very large
eigenvalues unlike the case of the algorithm in \cite{dicksonpal}
where the linearization introduces unwanted defective eigenvalues at
infinity which may create numerical problems (although such problems can be effectively amended by 
a structured refinement \cite{nuovo}).

If $nk$ is even, then by means of simple formal manipulations one may show that
$q(z):=x(z)^{-nk/2}\cdot p(x(z))$  is a polynomial in $z$,
where $x(z)=(z+\sqrt{z^2-4})/2$ or $x(z)=(z-\sqrt{z^2-4})/2$, i.e.,
$z(x)$ is one of the two branches of the inverse function of $z(x)=x+1/x$.

Moreover, taking derivatives in the latter equation leads to an
explicit expression for the Newton correction $q(z)/q'(z)$ given in
terms of $p(x)/p'(x)$
\begin{equation}\label{magicpal}
\frac{q(z)}{q'(z)}=\frac{1-1/x^2}{p'(x)/p(x) - nk/(2x)},
~~~p'(x)/p(x)=\hbox{tr}(P(x)^{-1}P'(x)).
\end{equation}

This equation enables one to apply the EAI to the polynomial
$q(z)$. Once its roots $z_1,\ldots,z_{nk/2}$ have been computed, the
eigenvalues of $P(x)$ are given by the pairs $(x_i,1/x_i)$ which are
the roots of the quadratic polynomial $x^2-z_ix+1$. This approach has
the advantage to work with an approximation vector of half the size
and to deliver the solution as pairs $(x,1/x)$.

It is important to point out that the applications $z\to
x=(z\pm\sqrt{z^2-4})/2$ is ill conditioned at $z=\pm 2$. Therefore,
loss of accuracy is expected near $x=\pm 1$. In this case, a
refinement step is advisable. Such a refinement may be implemented by
an unstructured version of the EAI, by the naive strucutred EAI, or with other structured refinement
methods \cite{nuovo}.

If $nk$ is odd, then $-1$ is necessarily an eigenvalue of the
palindromic PEP and there is no need to approximate it. To calculate
approximations of the remaining $nk-1$ eigenvalues, there are two
possible strategies.

As a first possibility, one may consider the new matrix polynomial
$Q(x)=(x+1)P(x)$ which has even degree. The eigenvalues of $Q(x)$ are
those of $P(x)$ and the eigenvalue $-1$ with multiplicity increased by 
$n$. Therefore, the previously described technique can be
applied.  Only $nk-1$ roots of $\det(Q(x))$ are needed,
because $n+1$ roots are a priori known to be equal to $-1$. Thus, one could 
apply the EAI \eqref{eaj} or \eqref{eags} 
with an approximation vector $y$ of $n(k+1)$ components of which $n+1$
are set equal to $-1$ in order to immediately achieve implicit
deflation of the roots; or, working in the variable $z$ in order to extract the structure, 
the SEAI \eqref{magicpal} can be used setting $(n+1)/2$ starting points equal to $-2$.

A second possibility is to set $x:=w^2$ and to consider the
eigenvalues of the polynomial $Q(w)=P(x(w))$. The scalar polynomial
$q(w):=\det Q(w)$ has $2nk$ roots, which are the square roots of the
solutions of the original equation $p(x)=0$ that we have to solve. In
particular only $2nk-2$ roots are to be determined, since $q(w)=0$ has
two known solution at $w=\pm i$. It is useful to set
$z:=(w+1/w)^2=x+1/x+2$. Defining
\[
\widetilde{q}(w):=\frac{q(w)}{w^{nk+1}+w^{nk-1}},
\]
 it is easy to check that $r(z):=\widetilde{q}(w(z))$ is a polynomial
 in $z$. Therefore we may restrict the attention to computing the
 roots of $r(z)$. Once they have been computed, the evaluation of the
 function $w(z)$ at these roots provides the roots of $q(w)$ . The
 evaluation of $x(w)$ at these latter roots yields the sought
 eigenvalues of $P(x)$.  In order to compute the roots of $r(z)$ we
 may apply the EAI to the polynomial $r(z)$. 
The following equations provides a tool to compute
 the Newton correction $r(z)/r'(z)$ needed by the EAI.
\[
\frac{r(z)}{r'(z)}=\frac{2 w
(1-1/w^4)}{ 
q'(w)/q(w) - [(nk+1)w^2+nk-1]/(w^3+w)},
\]
or in terms of the original variable $x$
\[
\frac{r(z)}{r'(z)}= \frac{1-1/x^2}{p'(x)/p(x)-[(nk+1)x+nk-1]/(2 x^2+2 x)}.
\]

At the moment we have no clear elements to say which of the two
possibilities is more convenient. We plan to investigate in this direction.

We conclude this subsection mentioning that also {\em antipalindromic}
and {\em anti-T-palindromic} polynomials ($\mathrm{Rev} P(x) = -P(x)$
and $\mathrm{Rev} P(x)=-P(x)^T$) have a $\{x,1/x\}$ symmetry. Their
determinants are pure palindromic if $n$ is even and antipalindromic
if $n$ is odd \cite{palsmith}.  The former case is exactly the same as
above. The latter case is also easy, because a scalar antipalindromic polynomial
is always equal to $x-1$ times a scalar pure palindromic polynomial. Moreover,
it is possible to prove  \cite{palsmith} that $1$ is always a root of a scalar
antipalindromic polynomial, and $-1$ is always a root of even-grade
antipalindromic polynomial, so according to the grade there are either one or
two exceptional eigenvalues with odd multiplicity. Therefore, it is
easy to extend our technique to this class.

\subsection{Hamiltonian/skew-Hamiltonian and even/odd}
An even (odd) polynomial is such that $P_j$ is symmetric for all even
(odd) values of $j$ and is skew-symmetric for all odd (even)
$j$. Similarly, the coefficients of a Hamiltonian/skew-Hamiltonian
polynomial are, alternatively, Hamiltonian and skew-Hamiltonian
matrices. The classes of even-dimensional even/odd polynomials are
easily mapped onto the classes of Hamiltonian/skew Hamiltonian
polynomials by a multiplication by $J$. Amongst the huge literature on
these classes of polynomials see, for instance, \cite{mehl, mewa, mewa2, urv} 
and the references therein.  Classical eigenvalue problems for Hamiltonian matrices \cite{vanloan} are a 
special case of skew-Hamiltonian PEPs.

The matrix polynomials belonging to these classes have eigenvalues
coming in pairs $(x,-x)$.
In particular, if $nk$ is odd, then either $x=0$ (if $P_0$ is skew-symmetric) or
$x=\infty$ (if $P_k$ is skew-symmetric)
is an eigenvalue. Notice that this is never the case for
Hamiltonian/skew-Hamiltonian polynomials, because they are only defined for even
$n$.

Let $z:=x^2$. Just like the T-palindromic case, also for even/odd
polynomials it is possible to follow the ideas exposed in
\cite{dicksonpal} and build a new matrix polynomial $M(z)$ whose
determinant is equal to $p(x(z)) \cdot p(x(z))$. The following result 
demonstrates the way it can be done for an even polynomial.

\begin{proposition}
Let $P(x)$ be an even matrix polynomial, and let $z=x^2$. Define
$B(z):=\frac{P(x(z))+P^T(x(z))}{2}$ and $C(z):=\frac{1}{x(z)}
\frac{P(x(z))-P^T(x(z))}{2}$, so that $P(x)=B(x^2)+xC(x^2)$. Then
$M(z):=\left(\begin{smallmatrix}
B(z) & z C(z)\\
C(z) & B(z)\\
\end{smallmatrix}\right)$ is a skew-Hamiltonian matrix polynomial such that $\det M(z) = [p(x(z))]^2.$ 
If $0 \neq x_0$ in $\mathbb{C}$ is 
an eigenvalue for $P(x)$ associated with a canonical set of 
Jordan chains of length $\ell_1, \dots, \ell_k$ then $x_0^2$ is an eigenvalue for $M(z)$ 
and its Jordan structure is the union of the Jordan structures of $P(x)$ at $x_0$ and at $-x_0$.

Moreover: 
\begin{enumerate}
\item concerning eigenvectors associated with any finite nonzero eigenvalue $x_0$, 
$P(x_0) v_0 = 0$ and $P(-x_0) w_0 = 0$ if and only if 
$v_1=(x_0 v_0^T, v_0^T)^T$ and $w_1=(-x_0 w_0^T,w_0^T)^T$ 
are two linearly independent eigenvectors such that $M(x_0^2) v_1 = 0 = M(x_0^2) w_1$;
\item $M(z)$ is a $2n \times 2n$ matrix polynomial of degree $\deg M = [(k+1)/2]$;
\item writing $M(z)=\sum_{j=0}^{\deg M} M_j z^j$, the relation $M_j=\left(\begin{smallmatrix}
P_{2j} & P_{2j-1}\\ P_{2j+1} & 
P_{2j}\end{smallmatrix}\right)$ holds for all $0 \leq j \leq \deg M$, where $P_j=0$ if $j<0$ or $j>k$;
\item if $k$ is odd, $M(z)$ has at least $n$ (respectively, $n+1$) eigenvalues at infinity if $n$ is even (respectively, odd).
          \end{enumerate}
\end{proposition}

The proof can be obtained by adapting the arguments used in
\cite{dicksonpal} for T-palindromic polynomials to the even case; we skip the details here.

Similar results can of course be obtained for odd and
Hamiltonian/skew-Hamiltonian matrix polynomials. Applying the EAI to
$M(z)$ allows us to extract the structure. An alternative approach,
that avoids possible issues about loss of accuracy for very large
eigenvalues (this time $M(z)$ has extra infinite eigenvalues only if $k$ is odd), is once again the implicit use of the squaring
transformation. Namely, if $nk$ is even (which is always satisfied for
Hamiltonian/skew-Hamiltonian polynomials), then defining $z:=x^2$ one
finds that $q(z):=p(x(z))$ is a polynomial for $x(z)=\sqrt z$ or
$x(z)=-\sqrt z$. Thus, $p'(x)/(2x p(x)) = q'(z)/q(z)$, so that the
Newton correction for the polynomial $q(z)$ is readily available
\[
q(z)/q'(z)=2xp(x)/p'(x)=2x/\hbox{tr}(P(x)^{-1}P'(x))
\]
 and the Ehrlich-Aberth algorithm can be implicitly applied to the
 polynomial $q(z)$ in order to compute its roots $z_1,\ldots,
 z_{nk/2}$. This way, the roots of $p(x)$ are readily available in pairs
 as $(\sqrt{z_i},-\sqrt{z_i})$.

In other situations, one eigenvalue is necessarily either $0$ (if
$P(x)$ is odd and $nk$ is odd) or $\infty$ (if $P(x)$ is even and $nk$ is odd); thus, there is no need to
approximate it. It may also happen that there is one uncoupled eigenvalue at $0$ and another one at $\infty$ 
(e.g. if $P(x)$ is odd, $n$ is odd and $k$ is even). In the case of an extra eigenvalue at $0$, to approximate the other eigenvalues 
one can notice that $q(z):=p(\sqrt z)/\sqrt z$ is a polynomial and that
$q'(z)/q(z) = (1/2x^2)(x p'(x)/p(x)-1)$. 
This yields the Newton correction for $q(z)$ as
\[
q(z)/q'(z)=2x/(p'(x)/p(x)-1/x)=2x/(\hbox{tr}(P(x)^{-1}P'(x))-1/x),~~~z=x^2,
\]
which enables one to apply the EAI to $q(z)$ by using an approximation
vector of length $(nk-1)/2$. As in the palindromic case, there is also the
alternative option to consider the polynomial $xP(x)$ which is even (odd) if $P(x)$
is odd (even). The new polynomial $xP(x)$ has $n$ additional eigenvalues at $0$ that
are known and can therefore be immediately deflated.

\subsection{Unified approach to any structure}

More in general, let $\C^*:=\C \cup \{\infty\}$ and let $f:\C^*
\rightarrow \C^*$ be any self-inverse function, that is $f(f(x))=x
\ \ \forall x \in \C^*$. An example is the subclass of rational
functions $f(x)=\frac{a x + b}{c x - a}$, which are self-inverse
whenever $a^2+bc \neq 0$. If we additionally require $f$ to be
analytic, having such a form is not only a sufficient condition, but
it is also necessary (unless $f(x)=x$) for $f$ to be
self-inverse. This follows from the fact that M\"{o}bius functions
(i.e., rational functions of degree $1$) are the only automorphisms of
$\C^*$.

Suppose that, because of some structure in the coefficients of $P(x)$,
all eigenvalues come in pairs $\{\lambda,f(\lambda)\}$. Eigenvalues
such that $\lambda=f(\lambda)$ are called \emph{exceptional}, and are
allowed to appear with any multiplicity. 

Such a possibility justifies the requirements that 
\begin{enumerate}
\item $f(x)$ is analytic, so that either it is the identity function or it 
has a finite number of fixed points, and
\item there is a way to identify which exceptional eigenvalues, if any,
appear with odd multiplicity.
\end{enumerate}
In fact, exceptional eigenvalues can otherwise
become a problem. For instance, consider a
real matrix polynomial associated to the non-analytic function $f(x)=x^*$: the method meets problems
in this case. The reason is that, since all the real line is exceptional, there
is no way to state \emph{a priori} which exceptional (i.e. real)
eigenvalues, if any, appear without being part of a complex conjugate
eigencouple.

If $f$ is analytic, the implicit change of variable method that we
have described for the special cases $f(x)=-x$ and $f(x)=1/x$ can be
generalised in the following way.  Suppose first that $a\neq 0$, and
define $z(x):=\frac{a x^2+b x}{cx-a}=xf(x)$. Let $x(z)$ denote any 
of the two branches of the inverse function of $z(x)$. Then if there are no eigenvalues with odd multiplicity one
can see that $q(z):=\frac{p(x(z))}{(cx(z)-a)^{nk/2}}$ is a polynomial;
therefore, the EAI can be applied to $q(z)$ and the eigenvalues can be
found inverting the rational function $z(x)$. Otherwise (e.g. if $nk$ is odd) there
must be some exceptional eigenvalues that can be treated with
techniques akin to those described for the special cases previously
considered. If on the contrary $a=0$, let $z(x):=\frac{c
  x^2+b}{cx}=x+f(x)$. Once again if all the eigenvalues come out in couples then
$q(z):=\frac{p(x(z))}{x(z)^{nk/2}}$ is a polynomial; otherwise one can simply deal with the exceptional eigenvalues with known 
odd multiplicity
by using techniques analogous to those described in the previous
subsections. Notice that the fixed points of $f(x)$ may lead to computational problems, since they are double roots in the equation
$z(x)=\zeta$. Refinements of some kind are advisable there. See also
\cite{nuovo}.

The explicit method may also be extended to the general case. This is the subject of a future research project.

\section{Numerical experiments}\label{numex}
We have performed extensive numerical experiments in order to check the efficiency and the accuracy of our implementation of the EAI. Further tests have been performed to confirm the ability of our method to exploit structures in the coefficients and to respect structures in the spectrum when approximating it. 
\subsection{Efficiency}
The complexity of the proposed algorithm is of order $tn^3+tkn^2$, where $t$ is the total number of times that a trace computation is needed before (vectorial) convergence. There is empirical evidence that $t$ heavily depends on the choice of the initial approximation. For the case of scalar polynomials, the use of suitable strategies \cite{B, BF} leads to a linear dependence of $t$ with respect to the total number of roots.

Experiments we made with our implementation, with starting points determined by the Newton polygonal, suggest that this is also the case of the EAI applied to a matrix polynomial. This means that the computational complexity of the EAI is $\mathcal{O}(kn^4+k^2n^3)$, leading to great computational advantages for $k \gg \sqrt{n}$. As noticed in \ref{linear}, if on the contrary $k \lsim \sqrt{n}$ other more focused implementation of the EAI are possible, with cubic efficiency in $kn$.

In order to confirm such predictions, we have compared our implementation of the EAI and Matlab\footnote{Matlab is a registered trademark of The MathWorks, Inc.}'s QZ implementation \emph{polyeig} on random matrix polynomials of high degree and small dimension. The experiments have been performed on a machine with CPU Intel Xeon 2.80GHz and system Linux Debian 6.02. 
 For very large values of $k$, we did not actually run \emph{polyeig} due to the very large forecast computation times, but we extrapolated the times from the other experiments; in fact, when doubling the value of $k$ we can expect that the running time of the QZ algorithm grows approximately with a factor $8$. Such extrapolated values are marked with a $^*$ in the following tables.
\medskip

\begin{tabular}{|c|c|c|}
\hline
\multicolumn{3}{|c|}{Computation times for $n=2$}\\
\hline
k & EAI & polyeig\\
\hline
50 & 0.018 s & 0.015 s\\
\hline
100 & 0.044 s & 0.064 s\\
\hline
200 & 0.111 s & 0.369 s\\
\hline
400 & 0.360 s & 4.35 s\\
\hline
800 & 1.29 s & 51.9 s\\
\hline
1600 & 4.76 s & 437 s\\
\hline
3200 & 18.4 s & $\mathcal{O}$(50 min)$^*$\\
\hline
\end{tabular} \begin{tabular}{|c|c|c|}
\hline
\multicolumn{3}{|c|}{Computation times for $n=5$}\\
\hline
k & EAI & polyeig\\
\hline
20 & 0.062 s & 0.010 s\\
\hline
40 & 0.121 s & 0.057 s\\
\hline
80 & 0.312 s & 0.370 s\\
\hline
160 & 0.920 s & 4.39 s\\
\hline
320 & 2.92 s & 44.0 s\\
\hline
640 & 10.3 s & 398 s\\
\hline
1280 & 38.1 s & $\mathcal{O}$(50 min)$^*$\\
\hline
2560 & 148 s & $\mathcal{O}$(7 hours)$^*$\\
\hline
5120 & 575 s & $\mathcal{O}$(2 days)$^*$\\
\hline
\end{tabular}
\medskip

The values in the tables above are in agreement with our prediction that the computation time should asymptotically grow as $k^2$. Moreover, the experimentation also confirms that, for a given value of $n$, the ratio of the time needed by EAI with respect to the time needed by the QZ algorithm exhibits an asymptotic growth that is approximately linear in $k$. This effect is taken to the extreme in the case $n=5$, $k=5120$. Had we used Matlab's \emph{polyeig}, it would have taken several days of computation time on our machine to solve such a problem. Our implementation of the EAI gave the approximated eigenvalues in less than $10$ minutes.
\subsection{Accuracy}
In order to test the accuracy of our implementation we used the Matlab toolbox NLEVP\cite{nlevp}. This toolbox has been recently proposed by its authors as an interesting set of benchmark problems that may be used as a standard test for new methods for nonlinear eigenvalue problems. It contains data coming from practical applications as well as model problems known to have peculiar properties.

Amongst the many nonlinear eigenproblems contained in NLEVP, we have selected all the square polynomial eigenproblems with $n < 25k^2$. We discarded PEPs with a larger ratio $n/k^2$ because they could be better dealt with by a different implentation of EAI, via a preliminary linearization. The test suite selected with this criterion consists of $29$ problems plus the $2$ problems butterfly and wiresaw1 that, being structured, will be treated in the next subsection.

In all the parameter-dependent problems in the library the default values of the parameters were selected. All methods were directly applied 
to the original matrices as saved in the library, without preprocessing them with 
any scaling. 
 Forward errors are evaluated by comparing the approximations with either theoretically known values, when available, or values computed in variable precision arithmetic (VPA) with Matlab's symbolic toolbox.

The graphs below are in logarithmic scale. Whenever the absolute error for a certain eigenvalue $\lambda$ appeared to be numerically zero, i.e. it was less than $\lambda$ times the machine epsilon $\epsilon =2^{-52} \simeq {\tt 2.22e-16}$, we formally set it equal to $\frac{\lambda \epsilon}{2}$. Only absolute errors for the finite eigenvalues are shown in the figures.

For the $3$ problems with $k\geq 3$, the eigenvalue forward errors
where computed for both the EAI and the QZ method (as implemented in
\emph{polyeig}).  Absolute errors for our implementation of the EAI
are marked with a red $*$ symbol, while absolute errors for
\emph{polyeig} are marked with a blue $+$ symbol. 
For this set of experiments, we
picked starting points on the unit circle. In our experience the order
of magnitude of the forward error is not significantly affected by the
choice of the starting points, even though for some problems other
choices led to slight improvements (not discussed here).

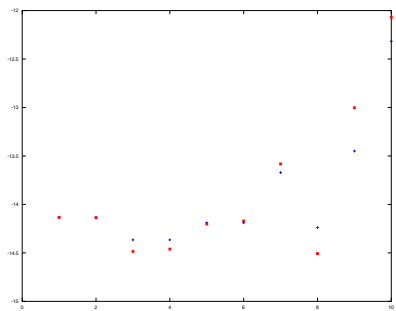
\begin{center}
\psone{orr_sommerfeld.epsc}{15}
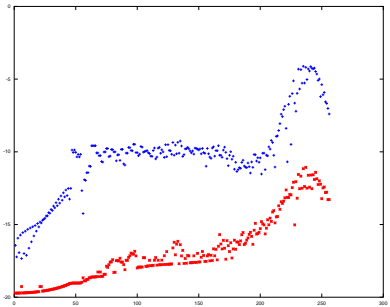
\fig{1}{Forward absolute errors for the problem orr sommerfeld}{0}
\end{center}
\begin{center}
\psone{plasma_drift.epsc}{15}
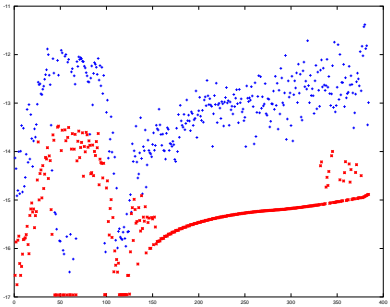
\fig{2}{Forward absolute errors for the problem plasma drift}{0}
\end{center}
\begin{center}
\psone{relative_pose_5pt.epsc}{15}
\fig{3}{Forward absolute errors for the problem relative pose 5pt}{0}
\end{center}

\medskip

For the $26$ problems with $k=2$, three methods were compared by computing their forward errors with the same method as above: 
\emph{polyeig} (blue $+$ symbol), EAI (red $*$ symbol) and the software \emph{quadeig} by Hammarling, 
Munro and Tisseur \cite{quadeig}, specifically designed for quadratic PEPs (black x symbol). Although we did not alter the 
coefficients given as input to any method, for  
most problems  \emph{quadeig} have performed scaling by the default settings of its internal algorithm, that 
prescribe scaling under certain conditions; see \cite{quadeig}. 

\medskip
\begin{center}
\pstwo{acoustic_wave_1d.epsc}{15}{acoustic_wave_2d.epsc}{15}
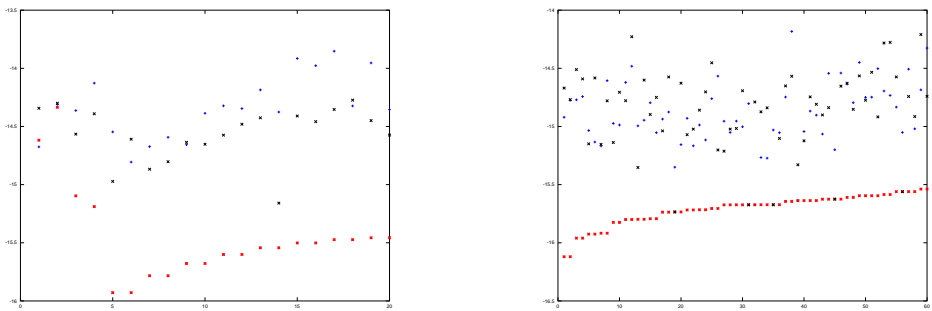
\fig{4}{Forward absolute errors for the problems acoustic wave 1d (left) and acoustic wave 2d (right)}{0}
\end{center}
\begin{center}
\pstwo{bicycle.epsc}{15}{bilby.epsc}{15}
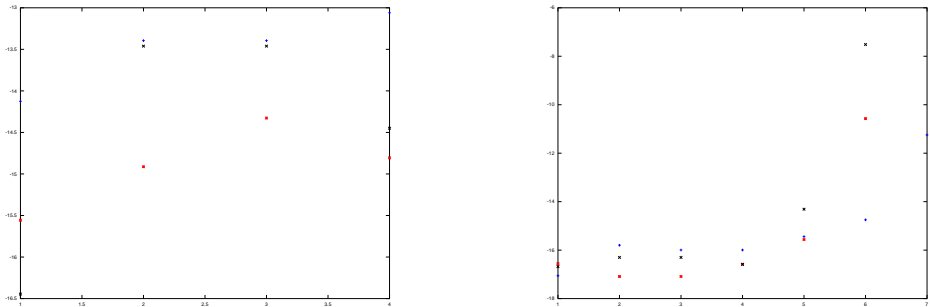
\fig{5}{Forward absolute errors for the problems bicycle (left) and bilby (right)}{0}
\end{center}
\begin{center}
\pstwo{cd_player.epsc}{15}{closed_loop.epsc}{15}
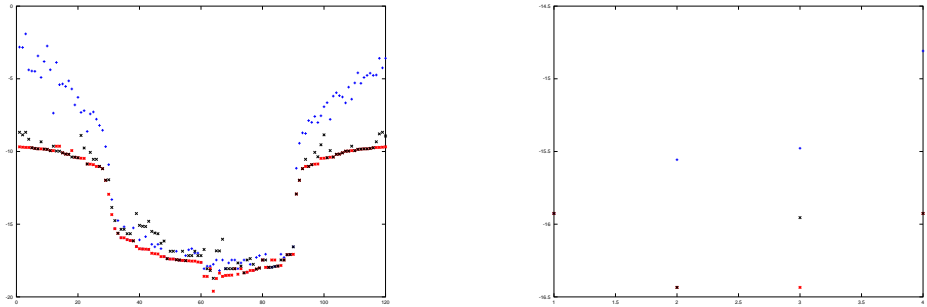
\fig{6}{Forward absolute errors for the problems cd player (left) and closed loop (right)}{0}
\end{center}
\begin{center}
\pstwo{dirac.epsc}{15}{gen_hyper2.epsc}{15}
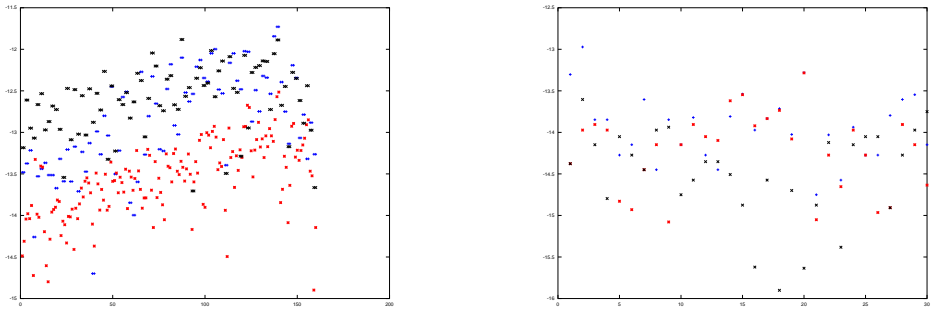
\fig{7}{Forward absolute errors for the problems dirac (left) and gen hyper 2 (right)}{0}
\end{center}
\begin{center}
\pstwo{hospital.epsc}{15}{intersection.epsc}{15}
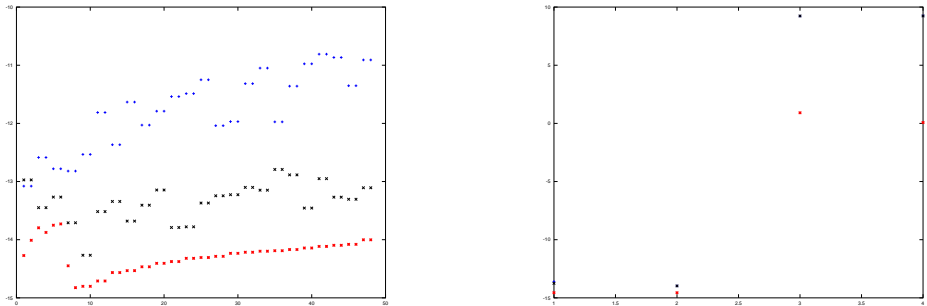
\fig{8}{Forward absolute errors for the problems hospital (left) and intersection (right)}{0}
\end{center}
\begin{center}
\pstwo{metal_strip.epsc}{15}{mobile_manipulator.epsc}{15}
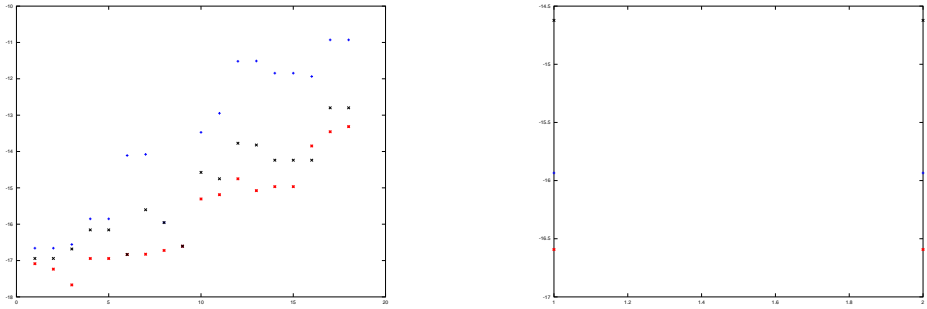
\fig{9}{Forward absolute errors for the problems metal strip (left) and mobile manipulator (right)}{0}
\end{center}
\begin{center}
\pstwo{omnicam1.epsc}{15}{omnicam2.epsc}{15}
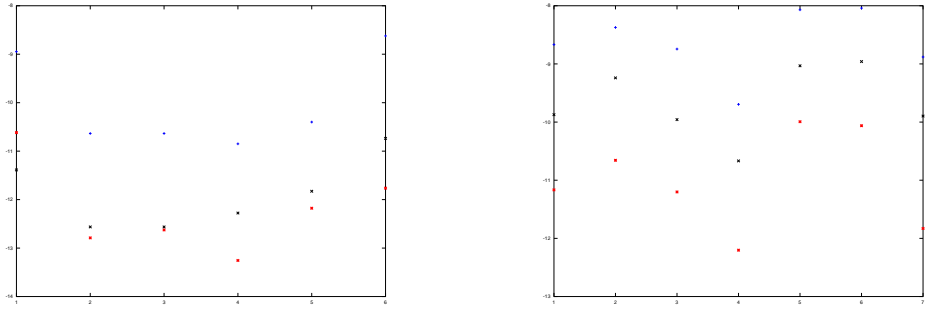
\fig{10}{Forward absolute errors for the problems omnicam1 (left) and omnicam2 (right)}{0}
\end{center}
\begin{center}
\pstwo{power_plant.epsc}{15}{qep1.epsc}{15}
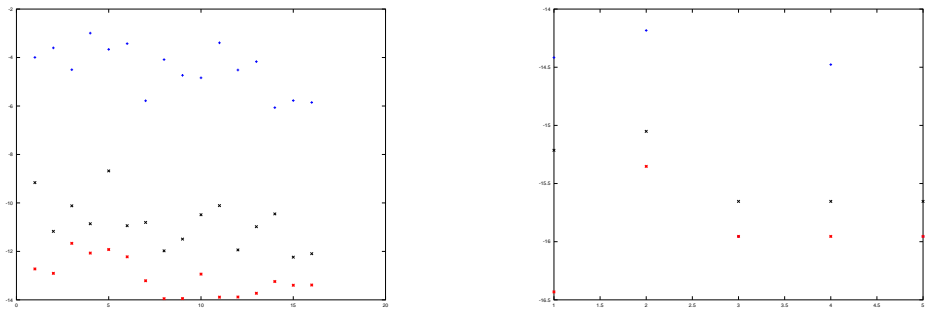
\fig{11}{Forward absolute errors for the problems power plant (left) and qep1 (right)}{0}
\end{center}
\begin{center}
\pstwo{qep2.epsc}{15}{qep3.epsc}{15}
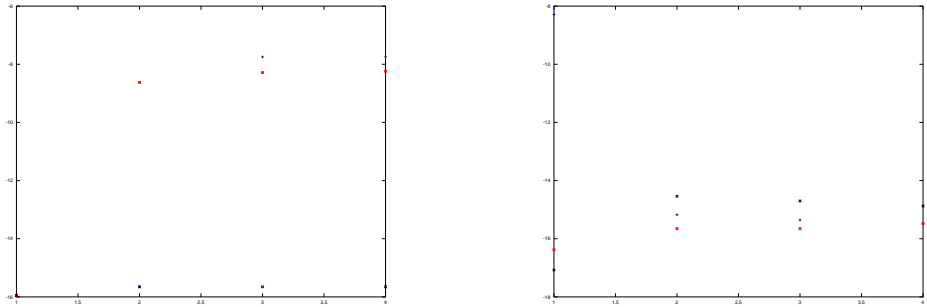
\fig{12}{Forward absolute errors for the problems qep2 (left) and qep3 (right)}{0}
\end{center}
\begin{center}
\pstwo{relative_pose_6pt.epsc}{15}{sign1.epsc}{15}
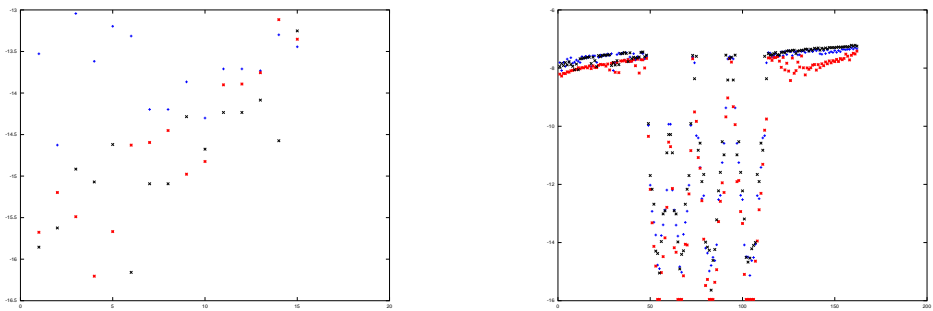
\fig{13}{Forward absolute errors for the problems relative pose 6pt (left) and sign1 (right)}{0}
\end{center}
\begin{center}
\pstwo{sign2.epsc}{15}{sleeper.epsc}{15}
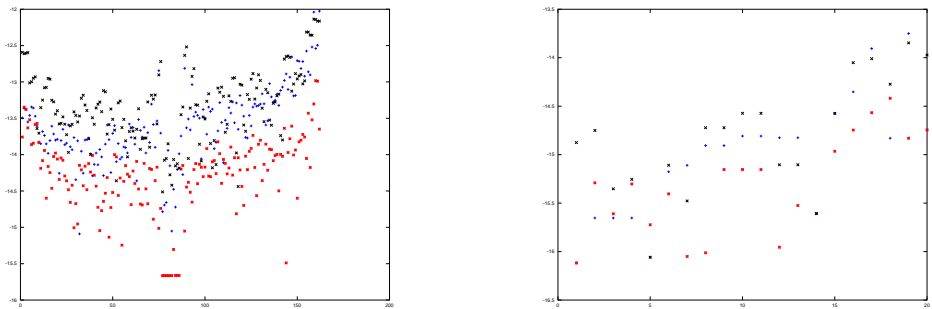
\fig{14}{Forward absolute errors for the problems sign2 (left) and sleeper (right)}{0}
\end{center}
\begin{center}
\pstwo{spring.epsc}{15}{spring_dashpot.epsc}{15}
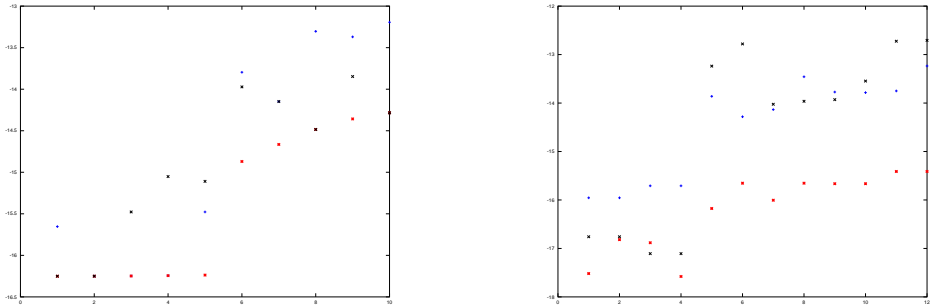
\fig{15}{Forward absolute errors for the problems spring (left) and spring dashpot (right)}{0}
\end{center}
\begin{center}
\pstwo{wing.epsc}{15}{wiresaw2.epsc}{15}
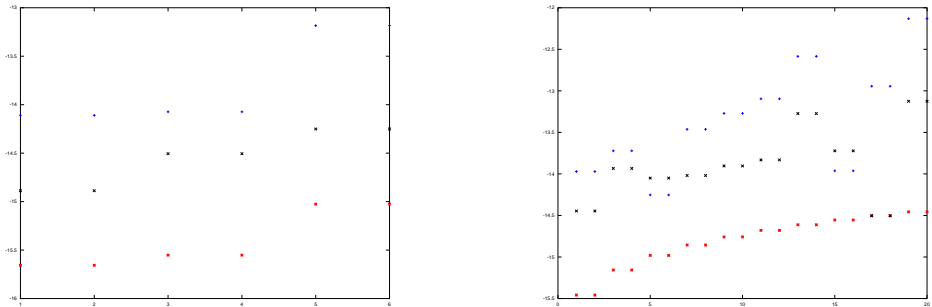
\fig{16}{Forward absolute errors for the problems wing (left) and wiresaw2 (right)}{0}
\end{center}

\medskip
As can be seen by the figures above, the approximations of the EAI are competitive, 
and often more accurate than the approximations of the QZ. In some cases, the improvement is remarkable. We report in the 
following table the maximal relative error and the average relative error for all the finite (i.e. neither numerically zero nor 
numerically infinite) eigenvalues of the $29$ considered problems, and for both the EAI and the QZ. The average relative error is 
defined as the geometric mean of all relative errors; numerically zero relative errors have been counted as relative errors equal 
to $\epsilon/2$. The values reported for the QZ for quadratic problems correspond to the 
algorithm, picked between \emph{polyeig} and \emph{quadeig}, that achieved the best performance in terms of average relative 
error for the given problem. As can be deduced by the above pictures and 
coherently with the results on backward errors presented in \cite{quadeig},
such best performance was achieved generally, but not always, by the latter algorithm.

\begin{tabular}{|c|c|c|c|c|}
\hline
Problem & \multicolumn{2}{|c|}{Rel. errors, EAI} & \multicolumn{2}{|c|}{Rel. errors, QZ}\\
\hline
 & Max. & Avg. & Max. & Avg.\\
\hline
acoustic wave 1d & 1.0e-14 & 2.1e-16 & 1.1e-14 & 1.7e-15\\
\hline
acoustic wave 2d & $\epsilon/2$ & $\epsilon/2$ & 3.2e-15 & 7.4e-16\\
\hline
bicycle & 1.0e-15 & 4.0e-16 & 7.6e-15 & 1.1e-15\\
\hline
bilby & 2.4e-14 & 3.5e-16 & 5.1e-15 & 1.8e-15\\
\hline
cd player & 5.3e-16 & 1.2e-16 & 4.0e-14 & 3.3e-16\\
\hline
closed loop & $\epsilon/2$ & $\epsilon/2$ & 3.4e-16 & 1.5e-16\\
\hline
dirac & 4.1e-14 & 5.9e-15 & 1.9e-13 & 2.9e-14\\
\hline
gen hyper 2 & 2.5e-14 & 9.3e-16 & 2.4e-15 & 4.9e-16\\
\hline
hospital & 2.7e-15 & 1.6e-16 & 2.0e-14 & 1.4e-15\\
\hline
intersection & 4.8 e-9 & 4.5 e-13 & 1.0 & 2.4e-8\\
\hline
metal strip & 6.3e-16 & 1.7e-16 & 2.3e-15 & 6.7e-16\\
\hline
mobile manipulator & $\epsilon/2$ & $\epsilon/2$ & 5.1e-16 & 5.1e-16\\
\hline
omnicam1 & 9.1e-11 & 1.2e-12 & 4.3e-9 & 6.4e-13\\
\hline
omnicam2 & 3.9e-10 & 2.3e-15 & 4.0e-9 & 1.2e-13\\
\hline
orr sommerfeld & 5.0e-12 & 9.1e-16 & 4.8e-5 & 2.8e-9\\
\hline
plasma drift & 3.4e-13 & 5.1e-16 & 1.3e-11 & 4.7e-14\\
\hline
power plant & 8.3e-14 & 1.1e-15 & 6.1e-11 & 1.9e-13\\
\hline
qep1 & 8.9e-16 & 1.7e-16 & 1.8e-15 & 5.1e-16\\
\hline
qep2 & 5.8e-9 & 5.3e-11 & 2.2e-16 & 1.9e-16\\
\hline
qep3 & 3.9e-9 & 1.0e-14 & 8.0e-10 & 3.2e-14\\
\hline
relative pose 5pt & 2.9e-14 & 6.8e-15 & 1.6e-14 & 6.3e-15\\
\hline
relative pose 6pt & 7.5e-14 & 1.9e-14 & 1.2e-13 & 1.4e-14\\
\hline
sign1 & 3.8e-8 & 1.1e-10 & 5.0e-8 & 3.9e-10\\
\hline
sign2 & 4.5e-14& 2.8e-15 & 3.1e-13 & 1.3e-14 \\
\hline
sleeper & 8.0e-16 & 2.8e-16 & 1.8e-15 & 5.6e-16\\
\hline
spring & $\epsilon/2$ & $\epsilon/2$ & 1.7e-15 & 3.5e-16\\
\hline
spring dashpot & 5.6e-15 & 3.1e-16 & 2.8e-13 & 1.4e-14\\
\hline
wing & $\epsilon/2$ & $\epsilon/2$ & 1.2e-15 & 8.1e-16\\
\hline
wiresaw2 & $\epsilon/2$ & $\epsilon/2$ & 2.4e-15 & 9.2e-16\\
\hline
\end{tabular}

\medskip

As the results above show, the EAI was generally able to improve the accuracy of the approximations with respect to the QZ method. The 
only problems where the EAI achieved an average performance worse than the QZ are QEP2 (loss accuracy on the multiple eigenvalue $1$ with 
respect to \emph{quadeig}; \emph{polyeig} has problems as well) and gen hyper 2.

The problems in NLEVP do not have high degree, so the condition $k^2 \gg n$ is not met. Therefore, in contrast with the high degree case, for those problems using the EAI as a primary algorithm does not bring advantages in term of computation time; on the contrary the implementation discussed in this paper is slower than QZ if $n \gg k^2$. Also, in these cases a linearization-based version of the EAI would be more efficient than the polynomial-based implementation discussed in the present paper. However, it is worth noticing that the numerical experiments showed that very often it improved the accuracy achieved by \emph{polyeig} and/or \emph{quadeig}; in many cases, the EAI is able to compute correctly all the digits of the eigenvalue up to machine precision. This suggests that, when $k^2 \lsim n$, it is possible to use the EAI as a refinement algorithm in order to improve the approximations obtained by the QZ. Using such values as starting points offers of course a very good choice, lowering the number of overall scalar iteration needed before convergence and therefore improving the efficiency of the EAI.
\subsection{Structured case}
Discussions and analyses of the behaviour of different structured versions of the EAI applied to structured PEPs have already appeared in \cite{dicksonpal} for palindromic polynomials and in \cite{paperino} for even/odd polynomials.

Here we further verify the reliability of our method with two tests. The first test relies on the NLEVP problems butterfly and wiresaw1, which are even \cite{nlevp}. The next figures compare four algorithms applied to these problems: \emph{polyeig}'s QZ (blue $+$ symbol), unstructured EAI (UEAI, red $*$ symbol), the structured matrix method URV applied to an even linearization \cite{urv} (black x symbol), and a structured version of the EAI relying on the change of variable $z=x^2$ method (SEAI, green o symbol).

\medskip
\begin{center}
\psone{butterfly.epsc}{20}
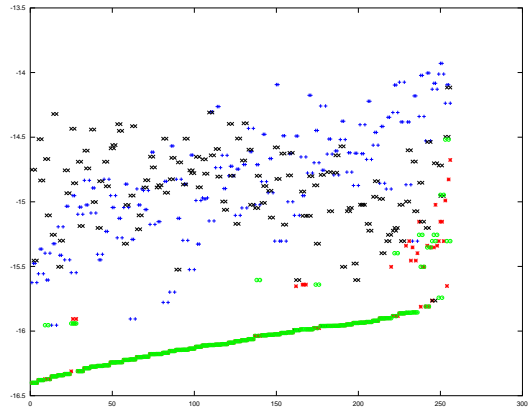
\fig{17}{Forward absolute errors for the problem butterfly}{0}
\end{center}
\begin{center}
\psone{wiresaw1.epsc}{20}
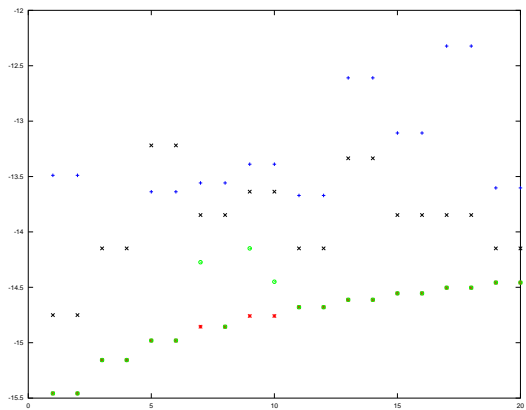
\fig{18}{Forward absolute errors for the problem wiresaw1}{0}
\end{center}
\medskip

It is clear from the figures above that for the problem butterfly the EAI was more accurate than the two matrix methods, but 
structured methods did not improve much the accuracy of each unstructured counterpart. 
This suggests that for these matrix polynomials the unstructured condition numbers for the eigenvalues are not much different from the structured condition numbers; therefore, structured methods do not improve much the accuracy, even though they improve the efficiency.

In the second test, a matrix polynomial $W(x)$ with $n=2$ and $k=10$ was built in such a way that its eigenvalues appear in couples of the form $\{\lambda,\frac{\lambda+1}{\lambda-1}\}$. In order to devise a problem not too easy to solve numerically, the determinant of the polynomial was designed to be Wilkinson-like: $\det(W(x))=const. \cdot \theta(x)\cdot \theta(\frac{x+1}{x-1}), \ \  \theta(x)=x\cdot\prod_{j=2}^{10}(x-j)$. The next figure shows the absolute errors of the computed approximations with respect to the known exact eigenvalues for three methods: QZ (polyeig, blue $+$ symbol), UEAI (red $*$ symbol) and SEAI relying on the change of variable $z=\frac{x^2+x}{x-1}$ (green o symbol). Numerically zero errors were formally set equal to $\epsilon/2$.
\medskip
\begin{center}
\psone{strange.epsc}{20}
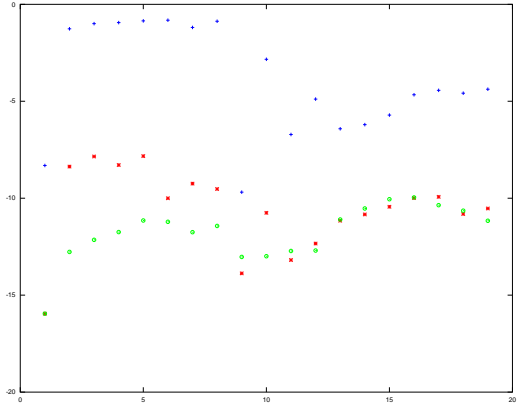
\fig{19}{Forward absolute errors for the structured problem $\det W(x)=0$}{0}
\end{center}
\medskip

The following table reports the relative errors of the three methods considered above for all the eigenvalues but $0$ (all the three 
algorithms detected the zero eigenvalue with an absolute error smaller than the machine epsilon). It also reports the relative error for
 the Matlab's function \emph{eig} (without scaling) applied to a suitable linearization, chosen according to the 
prescriptions of \cite{bestqz}. Notice that here all the nonzero eigenvalues have modulus $\geq 1$, so the suggested (near-to-optimal)
 linearization in the space $\mathbb{DL}$ according
to \cite{bestqz} would in principle be the pencil in $\mathbb{DL}$ corresponding to the ansatz vector $e_1$ (see \cite{vecspa} 
for further details). Unluckily, since one
 eigenvalue is zero, that pencil is not a linearization at all \cite{vecspa}. Following the suggestions on conditioning of \cite{bestqz}
 and the theory on linearizations of \cite{vecspa}, we have therefore taken
 the slightly perturbed vector $e_1+ 2^{-23} e_{10}$ as an ansatz vector. The factor $2^{-23}$ has been heuristically chosen 
picking the integer $\alpha \leq 52$ that minimises the average relative error when applying \emph{eig} to the linearizations in $\mathbb{DL}$
associated to the ansatz vectors $e_1 + 2^{-\alpha} e_{10}$.

\medskip

\begin{tabular}{|c|c|c|c|c|}
\hline
Eigenvalue & R. e., polyeig & R. e., eig & R. e., UEAI & R. e., SEAI\\
\hline
-1 & 4.9e-9 & $\epsilon/2$ & $\epsilon/2$ & $\epsilon/2$\\
\hline
11/9 & 4.5e-2 & 4.8e-8 & 3.5e-9 & 1.4e-13\\
\hline
5/4 & 8.1e-2 & 2.2e-7 & 1.1e-8 & 5.7e-13\\
\hline
9/7 & 8.9e-2 & 4.2e-7 & 4.0e-9 & 1.4e-12\\
\hline
4/3 & 1.1e-3 & 4.1e-7 & 1.1e-8 & 4.6e-12\\
\hline
7/5 & 9.9e-2 & 2.2e-7 & 7.1e-11 & 5.1e-12\\
\hline
3/2 & 8.9e-2 & 6.3 e-8 & 3.8e-10 & 2.5e-12\\
\hline
5/3 & 3.9e-2 & 8.9e-9 & 1.8e-10 & 1.1e-12\\
\hline
2 & 1.0e-10 & 4.9e-10 & 8.8e-12 &5.1e-14\\
\hline
2 & 7.3e-4 & 2.3e-12 & 6.7e-15 & 4.7e-14\\
\hline
3 & 4.3e-6 & 6.3e-12 & 2.2e-14 & 6.3e-14\\
\hline
3 & 6.5e-8 & 5.8e-12 & 1.5e-13 & 6.7e-14\\
\hline
4 & 9.6e-8 & 1.5e-10 & 1.7e-12 & 2.0e-12\\
\hline
5 & 1.3e-7 & 8.1e-10 & 2.9e-12 & 5.9e-12\\
\hline
6 & 3.2e-7 & 2.2e-9 & 6.1e-12 & 1.5e-11\\
\hline
7 & 3.1e-6 & 3.3e-9 & 1.5 e-11 & 1.6e-11\\
\hline
8 & 4.6e-6 & 2.7e-9 & 1.5e-11 & 5.5e-12\\
\hline
9 & 2.9e-6 & 1.1e-9 & 1.7e-12 & 2.5e-12\\
\hline
10 & 4.2e-6 & 1.7e-10 & 3.0e-12 & 6.9e-13\\
\hline
Average & 3.8e-5 & 9.6e-10 & 5.5e-12 & 4.1e-13\\
\hline
\end{tabular}

\medskip

We may conclude that on this structured problem the EAI outperforms the QZ method for what concerns accuracy. Apparently \emph{polyeig}
 struggles quite a bit here (which is coherent with the results of \cite{bestqz}), while the strategy 
of \cite{bestqz} works better, but still worse than the EAI. Although there are some approximations that do not benefit 
from the use of the structured version of the Ehrlich-Aberth algorithm, the SEAI has an overall advantage in accuracy over the UEAI, besides the obvious efficiency advantage. The computation time for the SEAI was about one third of the computation time for the UEAI.

\section{Conclusions and lines of future research}
We have proposed and tested a generalisation of the Ehrlich-Aberth
method to polynomial eigenvalue problems. Both theoretical arguments
and numerical experiments show that the Ehrlich-Aberth algorithm is
more efficient than customary method for high-degree PEPs, since its
complexity is only quadratic in $k$. Numerical experiments also
suggest that in many situations the new method provides more accurate
approximations, and therefore it may also be used as a refinement
method for low-degree PEPs.

On the other hand, problems arise in the treatment of multiple eigenvalues. Current research is focused on this issue. 
Currently, our algorithm exploits a heuristic device to deal 
with multiple eigenvalues at $0$ or $\infty$, which is the most common occurrence in practice. For the NLEVP problems, very good 
results were obtained also when multiple eigenvalues were present.

Another line of current research regards the possible use of other
root-finding algorithms. For instance, we can cite the modified Ehrlich-Aberth iteration \cite{macnamee, nourein}, the Durand-Kerner iteration 
\cite{macnamee, petkovic} or simultaneous root finders based on higher order methods in the Householder family, e.g. the Halley method 
\cite{petkovic}. It is advocated \cite{macnamee, petkovic} that some of the above mentioned method have order of convergence higher 
than $3$. However, in practice we did not see a significant improvement on the total number of scalar iterations with respect to the 
EAI; sometimes, performances were definitely worse than the EAI. Further work is needed to compare the various possibilities in special 
cases like structured matrix polynomials.

Important issues on which we plan to keep on working are improving the design of reliable stop conditions, a posteriori error bounds 
and choices of starting approximations. About the latter issue, in \cite{binnofsha} some results of \cite{sharify} on scalar polynomials 
will be generalised 
to matrix polynomials.

Finally, we have proposed and tested a structured version of the algorithm that is able to catch certain structures in the spectrum. 
New numerical experiments on the matter confirm and extend the results of \cite{dicksonpal, paperino} about the efficiency and 
the accuracy of the structured EAI.

\section{Acknowledgements}
The second author wishes to thank Federico Poloni who provided some help by delivering the eigenvalues computed by Matlab in VPA for some of the NLEVP problems.

\end{document}